\documentclass[12pt]{amsart}
\usepackage{amsmath,amsthm,epsfig,amssymb}
\usepackage{amsthm,amsfonts,amssymb,amscd}
\usepackage[all]{xy}
\newtheorem{thm}[subsection]{Theorem}
\newtheorem{prop}[subsection]{Proposition}

\newtheorem{lem}[subsection]{Lemma}
\newtheorem{corol}[subsection]{Corollary}
\newtheorem{rem}[subsection]{Remark}
\theoremstyle{definition}
\newtheorem{Def}[subsection]{Definition}

\newtheorem{exam}[subsection]{Example}
\newtheorem{proposition-definition}[subsection]{Proposition-Definition}

\newcommand{\ZZ}{{\mathbb Z}}

\newcommand{\PP}{{\mathbb P}}

\newcommand{\OOO}{{\mathcal O}}

\author{F. Laytimi}
\address{Math\'ematiques - b\^{a}t. M2, Universit\'e Lille 1,
F-59655 Villeneuve d'Ascq Cedex, France}
\email{fatima.laytimi@math.univ-lille1.fr}

\author{D.S. Nagaraj}

\address{Institute of Mathematical Sciences C.I.T. campus, Taramani, 
chennai 600113,India}
\email{dsn@imsc.res.in}

\title{Remarks on Ramanujam-Kawamata-Viehweg Vanishing Theorem}
\begin{document}
\begin{abstract} In this article we prove a general result on a nef vector bundle $E$ on a projective 
manifold $X$ of dimension $n$ depending on the vector space $H^{n,n} (X, E). $
It is also shown that $H^{n,n} (X, E)=0$ for an  indecomposable nef  rank 2 vector bundles $E$ on some specific type 
of $n$ dimensional projective manifold $X.$
The same vanishing shown to hold for indecomposable nef and big rank 2 vector bundles 
on any variety  with trivial canonical bundle.
\end{abstract}
 
\maketitle

\section{Introduction} \setcounter{page}{1}

Let  $X$  be a  smooth projective complex manifold of dimension $n.$ For any coherent sheaf $E$ on $X,$
we denote $H^{p,q} (X, E)$ the cohomology group $H^{q} (X, E
\otimes\Omega^p_X),$ where $\Omega^p_X$
is the sheaf of holomorphic  differential forms of degree $p$ on $X.$

  Akizuki-Kodaira-Nakano famous vanishing  theorem says:
 
    If $L$ is an ample line bundle on  a projective manifold $X$ of dimension $n,$ 
    then 
    
    $$H^{p,q} (X, L)=0  \ \ {\mbox {for}} \ \  p+q-n>0.$$ 
 The particular case $p=n$ is the  Kodaira vanishing theorem. The Kodaira vanishing theorem  was 
 extended to  nef and big   
 line bundle on a smooth surface by  Ramanujam \cite{Ra} and for higher dimension by Kawamata \cite{Ka} 
 and Viehweg \cite{Vi}.\\
 
  Ramanujam has given in \cite{Ra} an  example  showing that in general, one does not
expect  Akizuki-Kodaira-Nakano type vanishing result for nef and big line bundle. 
 
Le Potier \cite{Le} generalized the Akizuki-Kodaira-Nakano type vanishing  theorem to the case of ample
 vector bundle as follows: 
 
 If $E$ is an ample vector bundle of rank $r$ on a projective manifold $X$ of dimension $n,$ 
    then 
    \begin{equation}
     H^{p,q} (X, E)=0 \ \ \ { \mbox {for}} \ \ \ p+q-n> r-1.
    \end{equation}

The vanishing results of Ramanujam-Kawamata-Viehweg  and Le Potier  naturally led to ask the following question: 

Let $E$ be a nef and big vector bundle  of rank $r$ on a projective manifold $X$ of dimension $n.$ Is
\begin{equation}\label{eq0}
H^{n,q} (X, E)=0\ \ \ { \mbox {for}}\ \ \ q> r-1 ? 
\end{equation}
The example given by Ramanujam in \cite{Ra} shows that one can not expect in general 
"Akizuki-Kodaira-Nakano`` type of vanishing for nef and big line bundle. 
The same example can also be used to show that the question (\ref{eq0}) has a negative answer
(see Example (\ref{ex2})).

 Regarding the question (\ref{eq0}) for a  nef and big rank two vector bundle $E$ on a smooth surface $X$ 
 the only group which one hope  to vanish is
the group $H^{2,2} (X, E).$ In trying to investigate this problem  we obtained the following:
\begin{thm}\label{thm0}
 Let $E$ be a  nef vector bundle of rank $r$ on a projective manifold $X$ of dimension $n.$ 
Set $k(E):= {\rm dim}H^{n,n} (X, E).$ Then $k(E) \leq r$ and $E$ admits a trivial bundle
of rank $k(E)$ as quotient. In particular, $k(E)=r$ if and only if $E$ is isomorphic to trivial
vector bundle of rank $r.$
\end{thm}

\begin{corol}\label{cor0}
Let $E$ be an indecomposable nef vector bundle of rank $r$ on a projective manifold $X$ of dimension $n.$
Assume that $c_r(E) \neq 0.$ Then $H^{n,n}(X, E)=0.$
\end{corol}

For the case of rank 2 vector bundles we have the following:
\begin{thm}\label{thm1} 
Let $E$ be an indecomposable  nef vector bundle of rank $2$ on a projective manifold $X$ of dimension $n.$ 
If $H^{1} (X, \det (E))=0,  $ then $H^{n,n} (X, E)=0.$ 
\end{thm}

As a consequence  we obtain:
\begin{corol}\label{cor1} 
Let   $X$  be a Grassmannian  of dimension $n\geq 2$ 
or a complete intersection  of dimension $n\geq 3$ in a Grassmannian.\\ 
If  $E$ is an  indecomposable  nef vector bundle of rank $2$ on  $X, $  
then $H^{n,n} (X, E)=0.$ 
\end{corol}

\begin{corol}\label{cor2} 
Let $X$ be a projective manifold of dimension $n\geq 2$ with $K_X=\OOO_X. $ 
If  $E$ is an indecomposable  nef and big vector bundle of rank $2$ on  $X,$  
then $H^{n,n} (X, E)=0.$ 
\end{corol}

We recall a vanishing theorem of Schneider  \cite{Sc} related to nef and big vector bundle,
in a slightly different version  from the original one, but follows from the proof given there.

 \begin{thm}\label{Sc}
Let  $E$ (resp. $L$) be a vector bundle (resp. line bundle) on a projective manifold $X$ of dimension $n.$
If $E\otimes L$ is nef and big then 
 $$H^{n,q} (X,S^k (E)\otimes \det (E)\otimes L)=0,\ \ {\mbox {for}} \ \ q>0. $$
\end{thm}

\section{Notations and Definitions}

Throughout we work over the field of complex numbers. 

For a vector bundle $E$ on a projective manifold $X,$ we  will denote  by $E^{\vee}$ the dual of $E, $ 
$c_i(E) \in H^{2i}(X, \ZZ)$ is the $i$-th chern class of $E, $  
 $\PP(E)$ is the projective bundle whose  fiber over a point $x \in X$ is the projective space of $1-$dimensional quotients of 
the vector space $E_x, $ and $\OOO_{\PP(E)}(1)$ the universal quotient line bundle on $\PP(E).$

\begin{Def} Let $X$ be a projective
manifold of dimension $n.$ A line bundle $L$ on $X$ is called nef, if for every irreducible curve $C$ in $X $ 
degree of $L|_C$ is non negative. A nef line bundle $L$ is called big if $c_1(L)^n > 0.$

A vector bundle  $E$ on $X$ is said to be nef  if the line bundle $\OOO_{\PP(E)}(1)$ on $\PP(E)$ is nef. 
 
A nef vector bundle $E$ is said  to be big if 
$\OOO_{\PP(E)}(1)$ on $\PP(E)$ is big or equivalently 
$$s_n(E)= p_*(c_1(\OOO_{\PP(E)}(1))^{n+r-1}>0,$$ 
where $s_n(E)$ is the $n$-th Segre class of $E$ and $p :\PP(E) \to X$ be the natural projection. 
 \end{Def}

\section{Proof of the results}

First we recall some  results which we need.

\begin{prop}\label{prop1}\cite[Proposition 6.1.18 (i)]{LR1}
 A vector bundle  $E$ on $X$ is  nef if and only if the following condition is satisfied:\\
 Given any morphism $f: C \to X$ finite onto its image from an irreducible smooth curve 
 $C$ to $X,$ and given any quotient line bundle $L$ of $f^*(E),$ then one has ${\rm deg}L \geq 0.$
\end{prop}

\begin{lem}\cite[Proposition 1.16]{DPS}\label{lem1} 
Let $E$ be a nef vector bundle  on a projective manifold $X$ of dimension $n.$ 
If $\sigma$ is a non-zero section of  $E^{\vee}$     then $\sigma$ is nowhere vanishing on $X.$

\end{lem}

{\it Proof:} The proof given in \cite{DPS} uses analytic methods. Here we give an algebraic proof.
First we prove the lemma when $X$ is a curve. In this case if  $\sigma$ vanishes at some points, 
we get a positive degree line sub bundle of $E^{\vee}.$ 
By dualizing we see that $E$ has a line bundle quotient 
of negative degree. This is a contradiction to the Proposition(\ref{prop1}). 
Thus $\sigma$ is nowhere vanishing on $X.$

For the general case, assume   $\sigma$ vanishes at some points and dimension of 
$X$ is greater than one. Let $Z$ be  the
subscheme of $X$ defined by the vanishing
 of $\sigma$ and $I_Z$ denotes its sheaf of ideals.  The section $\sigma$ induces surjection
  \begin{equation}\label{eq1}
 \sigma :E \to I_Z \to 0. 
 \end{equation}
Let $C$ be a smooth curve in $X$ with the property $D= C \cap Z$ is a non-empty proper 
closed subscheme of $C.$ Then by restricting 
the surjective map $\sigma$ to $C$ and going modulo torsion we get a surjection: 
\begin{equation}\label{eq2}
 \tau :E|_C \to \OOO_C(-D) \to 0. 
 \end{equation}
Since $C$ is a smooth curve $\OOO_C(-D)$ is a line bundle of negative degree, which is
a contradiction to the fact that $E$ is nef. Hence we must have $Z = \emptyset.$  $\hfill{\Box}$

\begin{lem}\label{xi}\cite[see, Proposition 4.8.]{Xi}
If $E$ is a nef and big vector bundle on a K\"ahler manifold $X,$ 
then  the line bundle $\det (E)$ on $X$ is big.
\end{lem}

 The Dominance theorem [theorem 3.3] in \cite{LN} ensures that $\det (E)$ is nef.
 
We  also need to recall the proposition [Prop. 1.15 (iii)] in \cite{DPS}. 
 We will state it in a different version, which follows immediatly from the proof given there.
 \begin{lem}\label{nice}
 Let 
 $$0\to F\to E\to Q\to 0$$ 
 be an exact sequence of  holomorphic vector bundles and  
 $rank(E)=r, rank(F) =f.$\\  
 If $\wedge^{r-f+1}E \otimes \det Q^{-1}$ is nef  
 (resp. ample), then $F$ is nef (resp. ample).  
 \end{lem}

\bigskip

 {\bf Proof of Theorem(\ref{thm0})}:

\bigskip
The proof is by induction on the $rank(E)=r.$
If  $r=1$  and $k(E)=0$ then there is nothing to prove. \\
If $k(E) > 0,$ then by Lemma(\ref{lem1})
there is a non zero homomorphism 
$$ \sigma : \mathcal{O}_X \to E^{\vee} $$
which is nowhere
vanishing. This implies that $E$ is a trivial bundle of rank one. Since $k(\OOO_X) =1,$
the Theorem follows in this case.

Let $r>1.$ We assume our Theorem holds for all nef vector bundles of rank less than or equal $r-1.$
Again,  if $k(E)=0$ there is nothing to prove. So we assume $k(E) > 0.$ Then applying Lemma(\ref{lem1})
we get an exact sequence 
\begin{equation}\label{eq4}
  0 \to \OOO_X \to E^{\vee} \to F^{\vee} \to 0, 
\end{equation}
where $F^{\vee}$ is a dual of vector bundle $F$ of rank $r-1.$
Dualizing (\ref{eq4}) we get an exact sequence
\begin{equation}\label{eq5}
 0 \to F \to E \to \OOO_X \to 0.
\end{equation}
By Lemma(\ref{nice}) $F$ is a nef vector bundle. Now since  $rank(F)=r-1,$  
we have by induction assumption $k(F)\leq r-1$ and $F$ admits a trivial quotient of rank $k(F).$ 
This implies by duality
$F^{\vee}$ admits trivial subbundle of rank $k(F).$ 
Now from the long cohomology exact sequence 
\begin{equation}\label{eq6}
  0 \to H^0(X,\OOO_X) \to H^0(X, E^{\vee}) \to H^0(X, F^{\vee}) \to \cdots. 
\end{equation}
assosiated 
to the exact sequence (\ref{eq4}), we deduce 
 $$k(E)-1\leq k(F)\leq r-1,  $$
and the image of global sections of $E^{\vee}$ generate a trivial subbundle 
$V$ of rank $k(E)-1$ in $F^{\vee}.$
Taking the inverse image of this $V$   
we see 
that $E^{\vee}$ admits a subbundle $S^{\vee}$ of rank $k(E).$
Note that $S^{\vee}$ is an extension of ${\OOO_X}^{k(E)-1}$ by $\OOO_X.$ 
The dual $S$ of $S^{\vee}$ is nef, since it is an extension of trivial bundle of rank
$k(E)-1$ by a trivial bundle of rank $1.$
If $k(E)<r$ then it follows by induction $S$ is trivial. This proves the result.

If $k(E) = r$ then again by induction  $F ={\OOO_X}^{r-1}$ and all the sections of
$F^{\vee}$ lifts to sections of $E^{\vee},$   hence $E^{\vee}$ and $E$ are 
isomorphic to ${\OOO_X}^{r}.$ $\hfill{\Box}$
\bigskip

{\bf Proof of Theorem(\ref{thm1})}:

\bigskip

Assume $H^{n,n} (X, E)\neq 0,  $  then we get by Serre duality   $H^{0,0} (X, E^{\vee}) \neq 0.$
Let $\sigma$ be a non-zero section of 
$E^{\vee}.$ Since $E$ is nef by Lemma(\ref{lem1}) the section $\sigma$ is nowhere vanishing, 
and  gives an exact sequence
\begin{equation}\label{eq7}
 0 \to \OOO_X \to E^{\vee}\to {\det}(E)^{\vee} \to 0. 
 \end{equation} 
 This extension gives a class in the cohomology group $H^1(X,{\rm det}(E)).$ But by our assumption this group is
 zero and hence the extension splits. Thus $E^{\vee}$ splits and hence $E$ splits too, this  is a contadiction. $\hfill{\Box}$
 
 \bigskip

{\bf Proof of Corollary (\ref{cor0})}:
 
\bigskip
 If  $H^{n,n} (X, E)\neq 0, $ then by Theorem(\ref{thm0}) we get an exact sequence
 \begin{equation}
  0 \to F \to E \to \mathcal{O}_X \to 0.
 \end{equation}
This implies $c_r(E) = c_r(F) = 0,$ this is a contradiction. $\hfill{\Box}$

\bigskip

{\bf Proof of Corollary (\ref{cor1})}:
 
\bigskip

If $X$ is a Grassmannian of dimension $\geq 2$ or a complete intersection of dimension $\geq 3$ in  a Grassmannian,
 then for any line bundle $L$   $H^1(X,L) = 0.$ Hence if $E$ is an indecomposable vector bundle of rank two on $X,$ then
 the hypothesis of  Theorem(\ref{thm0}) holds for $E.$  $\hfill{\Box}$
\bigskip

{\bf Proof of Corollary (\ref{cor2})}

\bigskip

Assume $H^{n,n} (X, E) \neq 0.$ Since $E$ is nef and big, 
${\rm det} (E)$ is nef and big by the Lemma(\ref{xi}).
Hence we have an exact sequence:

\begin{equation}\label{eq8}
 0\to {\det(E)} \to E \to \mathcal{O}_X \to 0.
\end{equation}
 But  $K_X$ is trivial implies $H^1(X, \det(E))=0$ by Kawamata-Ramanujam-Viehweg vanishing theorem. Hence
that the exact sequence (\ref{eq8}) splits and hence $E$ is decompsable, which is a contradiction.
$\hfill{\Box}$

 \begin{rem} Corollary \ref{cor2} applies for example to 
 complex algebraic torus, $K3$  surfaces and Calabi-Yau  manifolds.
\end{rem} $\hfill{\Box}$

\section{Counter examples of Ramanujam}

\begin{exam}\label{ex1} The following example is due to Ramanujam \cite{Ra}.\\
Denote  $\PP^3$ blown up at a point by $X$ and
$\pi : X \to \PP^3$ be the natural morphism and $L = \pi^*(\OOO_{\PP^3}(1)).$
 Clearly the line bundle $L$ is  nef and big  and hence $H^1(X, \Omega^1_X\otimes L^{-1}) \neq 0.$
\end{exam}

\begin{exam}\label{ex2}
 Note that the variety  $X$ in the Example(\ref{ex1}) can be identified with $\PP(E)$ in such a way
that $L \simeq \OOO_{\PP(E)}(1),$ where $E=\OOO_{\PP^2} \oplus \OOO_{\PP^2}(1).$ Clearly the bundle $E$ on $\PP^2$ 
is nef and big and $H^{2,2}(\PP^2, E) \neq 0.$ This shows that one can not expect Le Potier type
vanishing result for nef and big vector bundle even for  $p=n.$ 

More general example: if $Y$ is a projective manifold of dimention $n$ and $H$ is an ample line bundle on $Y,$
then the vector bundle $E= \OOO_Y\oplus H$ is  nef and big vector bundle but  $H^{n,n}(Y,E)\neq 0.$
\end{exam}

\begin{rem} 
The non vanishing of $H^{1,1}(X, L^{-1})$ of Example(\ref{ex1}) 
can be deduced from the non vanishing   of the group $H^{2,2}(\PP^2, E)$ in Example(\ref{ex2}). 
Indeed:
$$H^2(X, \Omega^2_X\otimes L)\simeq H^{2,2}(\PP^2, \OOO_{\PP^2} \oplus \OOO_{\PP^2}(1))$$
by  Le Potier isomorphism  \cite[Lemma 8]{Le}.
\end{rem}

{\it Acknowledgements}: 
Last named author would like to thank University
of Lille1 and University of Artois at Lens and  would also like to 
thank  IRSES-Moduli program for their hospitality and the support.

\end{document}